\documentclass[11pt]{amsart}
\usepackage{amsmath, amscd}
\usepackage{amsfonts}
\usepackage{amssymb}
\usepackage{amsthm}
\usepackage{verbatim}
\usepackage[mathscr]{euscript}
\usepackage{graphicx}

\usepackage{xspace}
\usepackage{fullpage}

\def\ul{\underline}

\def\ds{\displaystyle}
\newcommand{\R}{\ensuremath{\mathbb{R}}\xspace}

\newcommand{\Pf}{{\bigskip\noindent\bf Proof.\quad} }

\newtheorem{theorem}{Theorem}[section] % numbered like the section
\newtheorem{lemma}[theorem]{Lemma} % numbered like the theorems
\newtheorem{proposition}[theorem]{Proposition}

\theoremstyle{definition} % styled differently... not italicized

\begin{document}
\title{A discrete model for nonlocal transport equations with fractional dissipation}
\author{Tam Do}
\email{Tam.Do@rice.edu}
\address{Dept. of Mathematics, Rice University, Houston, TX, USA 77005}
\date{\today}
%\date{}
\maketitle

\begin{abstract}
In this note, we propose a discrete model to study one-dimensional transport equations with non-local drift and supercritical dissipation. The inspiration for our model is the equation
$$
\theta_t + (H\theta) \theta_x +(-\Delta)^\alpha \theta =0
$$
where $H$ is the Hilbert transform. For our discrete model, we present blow-up results that are analogous to
the known results for the above equation. In addition, we will prove regularity for our discrete model which suggests supercritical regularity in the range $1/4<\alpha<1/2$ in the continuous setting.
\end{abstract}

\section{Introduction}

Before we introduce the equation we will primarily be concerned with, we will introduce the more complicated equation we hope our equation can model. Consider the equation
\begin{equation}
\label{eq:hild}
\theta_t (x,t)+ (H\theta)(x,t)\theta_x(x,t)+(-\Delta)^\alpha\theta (x,t)=0\,\quad (x,t)\in \R\times [0,\infty)
\end{equation}
where $H$ is the Hilbert transform and $(-\Delta)^\alpha$ where $\alpha>0$ is the fractional Laplacian. When $\alpha=0$, by convention, we take the dissipation to be absent in the equation. The equation can be viewed as a toy-model for the 2D surface quasi-geostrophic (SQG) equation
\begin{equation}
\label{eq:SQG}
\theta_t(x,t)+ u\cdot \nabla \theta(x,t) +(-\Delta)^\alpha\theta(x,t)=0,\quad u=\nabla^\perp(-\Delta)^{-1/2} \theta, \quad (x,t)\in \R^2\times [0,\infty)
\end{equation}
for which the question of regularity in the range $0\le \alpha<1/2$ is a longstanding open problem. See \cite{critsqg} and \cite{Kiselev} for recent progress and references.

\bigskip
By virtue of $\eqref{eq:hild}$ satisfying an $L^\infty$ maximum principle, we can roughly divide the analysis of $\eqref{eq:hild}$ into three categories according to the value of $\alpha$. For $\alpha>1/2$, the equation is $\textit{subcritical}$ in that one expects the dissipative term to dominate the nonlinearity and one has global regularity from functional analysis methods. For $\alpha=1/2$, the equation is $\textit{critical}$ in the sense that dissipative term and nonlinear term are thought to be equally balanced. For $0\le \alpha <1/2$, the equation is said to be $\textit{supercritical}$ in that the nonlinear term is thought to dominate the dynamics due to scaling considerations. It was shown in \cite{Blow2} and \cite{Blow1} that when $0\le \alpha<1/4$, solutions can develop finite time blow-up. In \cite{Dong1}, it was shown that for $\alpha\ge 1/2$ solutions exist for all time. In the supercritical range $1/4\le \alpha <1/2$, the question of global existence versus blow-up is open. It was shown in \cite{Do} that supercritical solutions derived in the vanishing viscosity limit, in a sense, become eventually regular. Also, it was shown that solutions supercritical by a logarithm are globally regular. We are not aware of any example where supercritical, with respect to the basic conservation laws, regularity has been proved for a fluid mechanics PDE.

\begin{figure}
\includegraphics[width=0.6\linewidth]{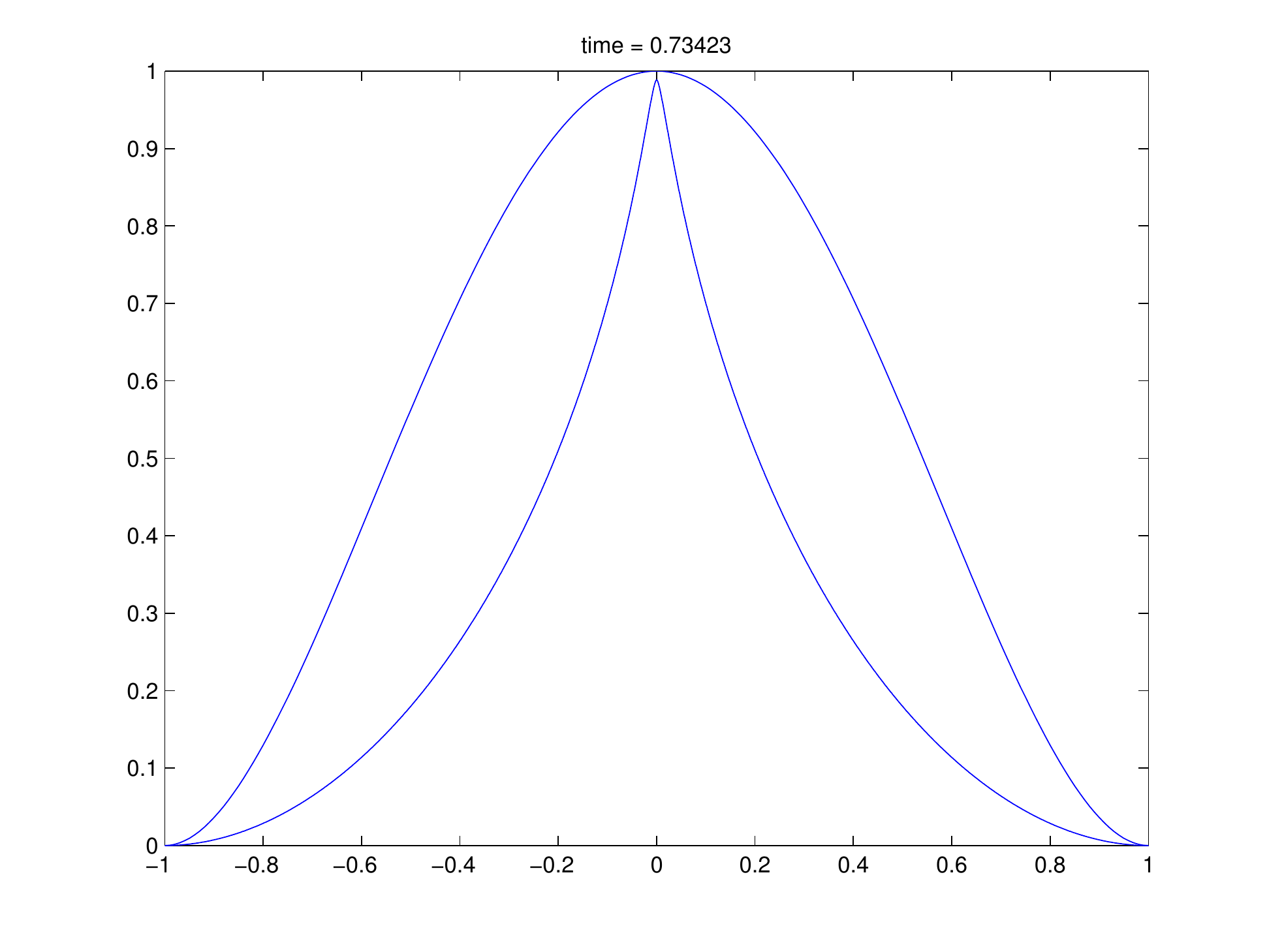}
\caption{Plot of the solution for the inviscid equation with initial data $\theta_0(x)=(1-x^2)^2\chi_{[-1,1]}(x)$}
\end{figure}

\subsection{Motivation for the model}

\bigskip
The proof of blow-up for $0\le \alpha <1/4$ relies on the following novel inequality that for $f\in C_c^\infty(\R^+)$ and $0<\delta<1$, there exists a constant $C_\delta$ such that
\begin{equation}
\label{eq:trick}
-\int_0^\infty \frac{f_x(x)(Hf)(x)}{x^{1+\delta}}\, dx \ge C_\delta \int_0^\infty \frac{f^2(x)}{x^{2+\delta}}\, dx.
\end{equation}
The proof given in \cite{Blow2} uses tools from complex analysis and is used to prove blow-up for even positive initial data with a maximum at $0$. In \cite{Kiselev}, another more elementary proof of blow-up for $0\le \alpha <1/4$ was given. The proof also goes by way of $\eqref{eq:trick}$ but without appealing to complex analysis. In particular, the complicated non-locality of the Hilbert transform is handled by the following key inequality from \cite{Kiselev}:

\begin{proposition}
 Suppose that the function $f(x)$ is $C^1$, even, $f'(x)\ge 0$ for $x>0$ and $f$ is bounded on $\R$ with $f(0)=0$. Then for $1<q<2$,
 $$
 Hf(x) \le \log(q-1)(f(qx)-f(q^{-1}x)).
 $$
\end{proposition}

\medskip
Consider initial data for $\eqref{eq:hild}$ that is $C^1$, even, and monotone decaying away from the origin. These properties are preserved by $\eqref{eq:hild}$. The proposition is then applied to $f(x)= \theta(0,t)-\theta(x,t)$. To derive our model, we will look at dyadic points. Consider the following system of ODEs
\begin{equation}
\label{eq:DE}
a_k'(t)=-(a_k-a_{k-1})^2(t)2^k, \quad a_k(0)=a_k^0, \quad k\ge 1.
\end{equation}
and set $a_0=0$.
This system serves as a discrete model for 
\begin{equation}
\label{eq:hil}
\theta_t (x,t)= -(H\theta)\theta_x(x,t).
\end{equation}
One can think of $a_k(t)$ as approximating $\theta(2^{-k},t)$. In particular, one could think that $-(a_k-a_{k-1})^2(t)2^k \approx \theta_x(2^{-k},t)$ and that, using Proposition 1.1 as inspiration, we shall take $(a_k-a_{k-1})$ as our approximation of $H\theta(2^{-k},t)$. The reason for considering a dyadic model is the nature of conjectured blow up for $\eqref{eq:hil}$. A cusp appears to form in finite time at the origin (see Figure 1) and it is our hope that studying $\eqref{eq:DE}$ will help us understand such phenomenon.

\bigskip
Now, we would like to formulate a discrete version of the fractional Laplacian $(-\Delta)^\alpha$. Recall for $\theta$ smooth enough,
$$
(-\Delta)^\alpha \theta(x)= P.V. \int_{-\infty}^\infty \frac{\theta(x)-\theta(y)}{|x-y|^{1+2\alpha}}\, dy ,
$$
\noindent
see \cite{Cordoba} for a derivation. If we take $\theta$ to be even then
$$
(-\Delta)^\alpha \theta(x)=\int_0^\infty \left[ \frac{1}{|x-y|^{1+2\alpha}}+ \frac{1}{|x+y|^{1+2\alpha}}\right] (\theta(x)-\theta(y))\, dy
$$
Let $x=2^{-k}$ and $y=2^{-n}$. If $k<n$ then $\ds  \frac{1}{|x-y|^{1+2\alpha}}+ \frac{1}{|x+y|^{1+2\alpha}} \approx 2^{(1+2\alpha)k}$. Similarly, if $k>n$, then $\ds  \frac{1}{|x-y|^{1+2\alpha}}+ \frac{1}{|x+y|^{1+2\alpha}} \approx 2^{(1+2\alpha)n}$. Also, $\theta(x)-\theta(y) \approx a_k-a_n$. Combining these observations, we can see that a reasonable model for fractional dissipation would be 
\begin{equation}
\label{eq:diss}
(\mathscr{L}^\alpha a)_k = \sum_{n=0}^{k-1} (a_k-a_n) 2^{2\alpha n}+ \sum_{n=k+1}^\infty (a_k-a_n) 2^{2\alpha k} 2^{k-n}.
\end{equation}
In our discrete formulation, we have ignored the tails of the integral defining $(-\Delta)^\alpha$. Adding the tails into (6) (i.e. letting the first sum go to $-\infty$) doesn't change the main results of this note and just adds more complications to the calculations.

\bigskip
Thus, our model for $\eqref{eq:hild}$ is
\begin{equation}
\label{eq:DDE}
a_k'= -(a_k-a_{k-1})^2 2^k - (\mathscr{L}^\alpha a)_k.
\end{equation}
$$
a_0'=-(\mathscr{L}^\alpha a)_0
$$
For simplicity of notation, we will sometimes omit the $\alpha$ and simply write $\mathscr{L}$ instead.
It should be noted that our model can also serve as discrete model for other non-local transport equations with fractional diffusion such as 
$$
\theta_t(x,t)= u(x,t)\theta_x(x,t)- (-\Delta)^\alpha (x,t)
$$
where
$$
u(x,t)= \left\lbrace \begin{array}{ll} \theta(x,t)-\theta(2x,t) & x\ge 0 \\
\theta(2x,t)-\theta(x,t) & x< 0 \end{array} \right.
$$
It is unknown whether solutions to this type of equations with general initial data can blow-up or exist globally in time in the supercritical range $0\le \alpha<1/2$. In \cite{Dylan}, a related kind of ``non-local'' Burgers type equation was studied and blow-up
was observed in the non-viscous case in certain situations.

\subsection{Outline of Main Results} In Section 2, we present several properties of solutions to $\eqref{eq:DDE}$ and show that solutions to $\eqref{eq:DE}$ blow-up, strengthening the analogy between the model and $\eqref{eq:hild}$. In Section 3, we prove an a-priori bound on solutions akin to a global in-time H\"older $1/2$ bound in the continuous setting. In Section 4, we use this bound to prove blow-up for $0<\alpha<1/4$ and global regularity for $1/4\le \alpha<1/2$, which can be considered the main result of this note. If this result were carried into the continuous setting, it would suggest that solutions to $\eqref{eq:hild}$ in the supercritical range $1/4\le \alpha<1/2$ are globally regular, contrary to natural scaling considerations.

\section{Local existence and properties of solutions}

\noindent
Define the space
$$
X^s=\{ \{a_k\}_{k=0}^\infty: \|a\|_X^s:= \sup_k |a_k| + \sup_{k\ge 1} |a_k-a_{k-1}|2^{sk} <\infty\}.
$$
\noindent
One could think of $X^s$ as being analogous to the H\"older spaces. However, the $X^s$ are made to deal with behavior near the origin.
Let $b_{k,s}= (a_k-a_{k-1})2^{sk}$.

\bigskip
\noindent
Before we show local existence for the full system, we will state some facts about $\mathscr{L}$ and its associated semigroup.

\begin{lemma}
\label{semigroup}
Let $0<\alpha<1/2$.

\bigskip
{\bf (a)} Suppose that for the index $k$, $b_{k,s}> c_s\|a\|_{X^s}$ where
$$
c_s = \frac{3}{4}(2^s-1)^{-1}\left(1 - \frac{1}{2^{s+1}-1}\right)<1.
$$
Then
$$
\left((\mathscr{L}^\alpha a)_k-(\mathscr{L}^\alpha a)_{k-1}\right)2^{sk} \ge C(\alpha)(2^{2\alpha k}-2^{2\alpha}) b_{k,s}
$$
where $C(\alpha)$ is a positive constant only dependent on $\alpha$.

\smallskip
	{\bf (b)}The operator $\mathscr{L}^\alpha$ generates a contracting semigroup $\ds  e^{-t\mathscr{L}^\alpha}$ on $X^s$ for all $s>0$.

\smallskip
	{\bf (c)} The following identity, which is analogous to $\int(-\Delta)^\alpha\theta=0$, is true
	$$
	\sum_{k=0}^\infty (\mathscr{L}^\alpha a)_k 2^{-k}=0.
	$$
\end{lemma}

\Pf
{\bf (a)} 
We will consider $k>1$ as the case for $k=1$ is similar. By direct computation,
\begin{eqnarray*}
(\mathscr{L} a)_k - (\mathscr{L} a)_{k-1} &=& \sum_{n=0}^{k-2} (a_k-a_{k-1})2^{2\alpha n} + 2(a_k-a_{k-1})2^{2\alpha (k-1)}  \\ \nonumber &\, & + 2^{(1+2\alpha)k} (1-2^{-(1+2\alpha)})\sum_{n=k+1}^\infty  (a_k-a_n)2^{-n} \\
&=& (a_k-a_{k-1}) \left( \frac{2^{2\alpha (k-1)}-1}{2^{2\alpha}-1}\right) + 2(a_k-a_{k-1})2^{2\alpha (k-1)}
 \\ &\,& + 2^{(1+2\alpha)k} \sum_{n=k+1}^\infty (a_k-a_n)(1-2^{-(1+2\alpha)}) 2^{-n} \\
& :=& I+ II +III
\end{eqnarray*}
By hypothesis, we have $c_sb_{j,s}< b_{k,s}$ for $j>k$. In terms of $a$, this means, for $n \ge k+1$,
$$
(a_n-a_k) = \sum_{j=k+1}^{n} (a_j -a_{j-1}) \le c_s^{-1}(a_k-a_{k-1}) \sum_{j=1}^{n-k} 2^{-js} = \frac{c_s^{-1}}{2^s-1}(a_k-a_{k-1}) (1-2^{-(n-k)s}).
$$
Then
\begin{eqnarray*}
III &\ge &- \frac{c_s^{-1}}{2^s-1}2^{(1+2\alpha)k} (a_k-a_{k-1})(1-2^{-(1+2\alpha)}) \sum_{n=k+1}^\infty (1-2^{-(n-k)s}) 2^{-n} \\
&=&  -\frac{4}{3} (1-2^{-(1+2\alpha)})(a_k-a_{k-1})2^{2\alpha k}.
\end{eqnarray*}
From this inequality and using that $0<\alpha<1/2$, it is easy to see that $II + III \ge 0$. Then
$$
\left((\mathscr{L} a)_k-(\mathscr{L}a)_{k-1}\right)2^{sk} \ge I \cdot 2^{sk} = C(\alpha)(2^{2\alpha k}-2^{2\alpha}) b_{k,s}
$$

\bigskip
{\bf (b)} Consider the system 
$$
a_k'=-(\mathscr{L} a)_k.
$$
\noindent
Using (a), it's not hard to see that $\ds \frac{d}{dt} \|a\|_{X^s }(t)< 0 $ and so $e^{-t\mathscr{L}}$ is a contracting semigroup.

\medskip
{\bf (c)} follows by explicit computation.

 $\Box $
 
\bigskip
\noindent
We will need the following version of Picard's theorem to prove local existence
\bigskip
\noindent
\begin{theorem} (Picard Fixed Point) Let $Y$ be a Banach space and let $\Gamma: Y\times Y\to Y$ be a bilinear operator
such that for all $a,b\in Y$,
$$
\| \Gamma(a,b)\|_Y \le \eta \|a\|_Y \|b\|_Y.
$$
Then for any $a^0\in Y$ with $4\eta \|a^0\|_Y < 1$, the equation $a=a^0+\Gamma(a,a)$ has a unique solution $a\in Y$ such that $\|a\|_Y \le 1/2\eta$.
\end{theorem}

\bigskip
\begin{theorem} (Local existence)
\label{existence}
Let $\{a_k(0)\}\in X^s$ where $s\ge 1$. Then there exists $T=T(\|a(0)\|_{X^s})$ such that there is a unique solution $\{a_k(t)\}\in C([0,T], X^s)$.
\end{theorem}

\Pf The argument is fairly standard and we will provide a sketch. By Lemma \ref{semigroup}, the operator $\mathscr{L}$ generates a contracting semigroup $\ds e^{-t\mathscr{L}}$ on $X^s$. Observe that a solution will satisfy
$$
a_k(t)= e^{-t\mathscr{L}}a_k^0 - \int_0^t e^{-(t-s)\mathscr{L}} \{2^k(a_k-a_{k-1})^2(s)\}\, ds
$$
In writing $\ds e^{-(t-s)\mathscr{L}}$ in the integral, we have slightly abused notation. The integrand is the $k$th element of $e^{-(t-s)\mathscr{L}}$ applied to the sequence $\{ 2^j(a_j-a_{j-1})^2(s)\}_{j=1}^\infty$.

Define a bilinear operator $\Gamma: X^s\times X^s \to X^s$ by
$$
\Gamma(a,b)_j(t)= -\int_0^t e^{-(t-s)\mathscr{L}} \{2^k(a_k-a_{k-1})(s)(b_k-b_{k-1})(s)\}\, ds
$$
\noindent
Define $\gamma(a,b)(s)$ by
$$
\gamma(a,b)_k (s)= -e^{-(t-s)\mathscr{L}}\{ 2^k(a_k-a_{k-1})(s)(b_k-b_{k-1})(s)\}
$$
Then after doing basic estimates, one has $\|\gamma(a,b)(s)\|_{X^s} \le C \|a(s)\|_{X^s} \|b(s)\|_{X^s}$ for $s\ge 1$. From this, we have an estimate on $\Gamma$ and choosing $t$ small enough, we can apply Picard.
$\Box$

\medskip
\noindent
Also, we have the following preservation properties:

\begin{lemma} Let $\{a_k(t)\}$ be a solution of $\eqref{eq:DDE}$ in $C([0,T], X^s)$, $s>1$. Suppose $\{a_k^0\}$ is non-decreasing in $k$ and is non-negative. 

\bigskip
		{\bf(a)} (Monotonicity and positivity) Then for all $t\le T$, $\{a_k(t)\}$ is non-decreasing in $k$ and is non-negative.
		
\bigskip
		{\bf{(b)}} (Max stays at $0$) For $t\in[0,T]$, we have that
$$
\sup_k a_k(t) = \lim_{k\to \infty} a_k(t)
$$

\bigskip
		{\bf(c)} ($\ell^\infty$ maximum principle) The supremum $\sup_k a_k(t)$ is non-increasing in $t$ and $a_0(t)$ is increasing in $t$. If $\{a_k(t)\}$ is a solution to $\eqref{eq:DE}$, then $\sup_k a_k(t)$ is constant in time \end{lemma}

\Pf 

{\bf(a)} Set $b_k(t)=(a_k-a_{k-1})(t)2^k$. For a contradiction, suppose there exists $j$ and a time $t_0$ such that $b_j(t_0)<0$. Choose a sufficiently small $\epsilon$ with $0<\epsilon<1$ such that $b_j(t_0)< -\epsilon/t_0$. Define a function $f$ by $\ds f(t)= \frac{-\epsilon}{2t_0-t}$ for $0\le t<2t_0$, so $b_j(t_0)<f(t_0)$.

\medskip
Let $g(t)=\inf_k b_k(t)$. Let $t_1$ be the first time $g$ crosses $f$ so $\inf_k b_k(t_1)= f(t_1)$. Because $b_k(t_1)\to 0$ as $k\to \infty$, the infimum is actually achieved for some index. Let $k_1$ be the first index for which $b_{k_1}(t_1)= f(t_1)$. By minimality of $t_1$, $\ds b_{k_1}'(t_1) \le f'(t_1)= -\frac{\epsilon}{(2t_0-t)^2}$. However, at time $t_1$, $b_{k_1}$ also satisfies
\begin{eqnarray*}
b_{k_1}' &=&-b_{k_1}^2+2b_{k_1-1}^2- \left((\mathscr{L} a)_{k_1}-(\mathscr{L}a)_{k_1-1}\right) 2^{k_1} \\
&\ge& -\frac{\epsilon^2}{(2t_0-t)^2} - \left((\mathscr{L} a)_{k_1}-(\mathscr{L}a)_{k_1-1}\right) 2^{k_1}
\end{eqnarray*}
By an argument analogous to Lemma \ref{semigroup}, the contribution from the dissipative terms on the right side of the inequality above 
is negative. From this, we arrive at a contradiction to the above inequality. Because $a_0$ is always non-negative, $a_k$ stays non-negative for all time by the monotonicity just proved.

\begin{comment}
%%% old wrong proof %%%%
Suppose, for a contradiction that the set $\{t>0: \mbox{for some k}\,\, a_k(t)\le a_{k-1}(t)\}$ is non-empty. The set has a minimum call it $t_0$. Let $k$ be the smallest index for which $a_k(t_0)=a_{k-1}(t_0)$. For the rest of the proof, everything will be evaluated at $t=t_0$. Then
$$
(a_k-a_{k-1})'=(a_{k-1}-a_{k-2})^2 2^{k-1} - \left( (\mathscr{L} a)_k - (\mathscr{L} a)_{k-1}\right).
$$
Since $a_{k-1} \ge a_{k-2}$, if we can show $(\mathscr{L} a)_k - (\mathscr{L} a)_{k-1} <0$, $(a_k-a_{k-1})' >0$ and we will have a contradiction. After some manipulations,
$$
(\mathscr{L} a)_k - (\mathscr{L} a)_{k-1} = 2^{(1+2\alpha)k} \sum_{n=k+1}^\infty (1-2^{-(1+2\alpha)})(a_k-a_n)2^{-n} <0
$$
since $0<\alpha<1/2$ and $a_k\le a_n$ for $n>k$. It is not hard too see that we cannot have $a_k=a_n$ for $n>k$ so the quantity in the above line is strictly less than $0$. 
\end{comment}

{\bf(b)} follows directly from {\bf(a)}.

\medskip
{\bf(c)} By {\bf(a)}, we see that $(\mathscr{L}a)_0(t)<0$ for all time $t\le T$ so $a_0(t)$ is increasing in $t$. Now, for $t\in[0,T]$
$$
\frac{d}{dt} \sup_{k>0} a_k(t) = \frac{d}{dt} \lim_{k\to \infty} a_k(t) =\lim_{k\to\infty} \frac{d}{dt} a_k(t) = \lim_{k\to\infty} \left(-(a_k-a_{k-1})^2(t) 2^k - (\mathscr{L}^\alpha a)_k(t)\right)
$$
Since we are on a compact time interval, the convergence above is uniform and so we are justified in interchanging limit and derivative. The limit of the first term on the far right is zero. Using that $\{a_k(t)\}\in X^s$ where $2\alpha <1 <s$,
$$
\lim_{k\to\infty}  (\mathscr{L}^\alpha a)_k(t) = \lim_{k\to \infty}\left(\sum_{n=0}^{k-1} (a_k-a_n)(t) 2^{2\alpha n}+ \sum_{n=k+1}^\infty (a_k-a_n)(t)2^{2\alpha k} 2^{k-n} \right) >0
$$ 
because the limit of the first term is positive by monotonicity and the limit of the second term is zero. This completes the proof.

$\Box $

\medskip
\noindent
In the spirit of \cite{Kiselev}, we show that solutions of $\eqref{eq:DE}$ blow up in finite time.

\bigskip
\begin{theorem} 
Let $\{a_k\}_{k=0}^\infty$ be a solution of $\eqref{eq:DE}$ in $X^s$, $s>1$ with initial data that is non-negative and non-decreasing in $k$. Then $\{a_k\}$ develops blow up in $X^s$ for every $s>1$ in finite time.
\end{theorem}
Define
$$
J(t)=\sum_{k=1}^\infty (a-a_k(t))2^{k\delta}
$$
where $0<\delta<1$. Observe that, $J\le \|a\|_{X^1} \sum_{k=1}^\infty 2^{k(\delta-1)}$. To prove we have blow up, we will show that $J$ must become infinite in finite time. We will need the following lemma

\noindent
\begin{lemma}
Let $0<\delta<1$. Suppose that $\{a_k\}\in X^s$, $s>1$, non-decreasing in $k$, and $a:= \lim_{k\to \infty} a_k <\infty$. Then
\begin{equation}
\label{ineq}
\sum_{k=1}^\infty (a_k(t)-a_{k-1}(t))^2 2^{k(\delta+1)} \ge C_0(\delta) \sum_{k=1}^\infty (a-a_k(t))^2 2^{k(\delta+1)}
\end{equation}
where $C_0(\delta)$ is constant depending only on $\delta$.
\end{lemma}

\Pf Choose $c>0$ such that $(c+1)^{-2}2^{\delta+1}>1$. We call $k$ "good" if $a_k-a_{k-1} \ge c(a-a_k)$ and bad otherwise.

\noindent
\medskip
{\bf Claim:} If the set of all good $k$ are finite, then both sides of $\eqref{ineq}$ are infinite.

\medskip
If the set of good $k$ is finite, then there exists $K$ such that for $k>K$, $a_k-a_{k-1} \le c(a-a_k)$ or equivalently $(a-a_k)\ge (c+1)^{-1}(a-a_{k-1})$. Then
$$
(a-a_k)^2 2^{k(\delta+1)} \ge (c+1)^{-2(k-K+1)} (a-a_K)^2 2^{k(\delta+1)} \to \infty
$$
as $k\to\infty$ by the choice of $c$, so the right side of $\eqref{ineq}$ diverges.

\medskip
Define $c_n$ by $a_k-a_{k-1}=c_k(a-a_k)$. Then $(a-a_{k-1})=(1+c_k)(a-a_k)$. Since $\lim_{k\to \infty}(a-a_k)=0$, $\prod_{k=K}^\infty (1-\frac{c_k}{c_k+1})= \prod_{k=K}^\infty \frac{1}{1+c_k}=0$, so $\sum_k c_k=\infty$. For the bad $k$,
$$
(a_k-a_{k-1})^2 2^{k(\delta+1)} = c_k^2 (a-a_k)^2 2^{k(\delta+1)}
$$
from which it is not hard to see that the left side of $\eqref{ineq}$ will be infinite as well. The claim is proven.

\medskip
Now suppose $k_{j-1}$, $k_j$ are good such that for $k_{j-1}<k<k_j$, $k$ is bad. Then $(a-a_{k-1})\le (1+c)(a-a_k)$, which implies
\begin{eqnarray*}
\sum_{k=k_{j-1}+1}^{k_j} (a-a_k)^2 2^{k(\delta+1)} &\le & \sum_{k=k_{j-1}+1}^{k_j} (1+c)^{2(k_j-k)} (a-a_{k_j})^2 2^{k(\delta+1)} \\
& = & (a-a_{k_j})^2 2^{k_j(\delta+1)} \sum_{k_{j-1}+1}^{k_j} (1+c)^{2(k_j-k)}2^{(k-k_j)(\delta+1)} \\
&=& (a-a_{k_j})^2 2^{k_j(\delta+1)}\sum_{n=0}^{k_j-k_{j-1}+1} (1+c)^{2n}2^{-n(\delta+1)} \\
& \le & C(\delta)  (a-a_{k_j})^2 2^{k_j(\delta+1)}
\end{eqnarray*}
where the last inequality comes from the choice of $c$. Treating all bad $k$ this way, the inequality $\eqref{ineq}$ follows. $\Box$

\bigskip
\noindent
{\bf Proof of Theorem 2.5}
Using the Lemma 2.4(c) and 2.6 as well as Holder's inequality, we have
\begin{eqnarray*}
\frac{d}{dt} J(t) &=& \sum_{k=1}^\infty (a_k(t)-a_{k-1}(t))^2 2^{k(\delta+1)} \ge C_0  \sum_{k=1}^\infty (a-a_k(t))^2 2^{k(\delta+1)} \\ & \ge &C_0 \left(\sum_{k=1}^\infty 2^{k(\delta-1)}\right)^{-1} \left(\sum_{k=1}^\infty (a-a_k(t))2^{k\delta}\right)^2 = C_1(\delta) J(t)^2
\end{eqnarray*}
\noindent
The result follows from Gronwall's inequality. $\Box $

\section{A-priori H\"older-$1/2$ bound}

The purpose of this section is to prove an a-priori bound for solutions of $\eqref{eq:DDE}$ from which the results of next section will follow. If this result of this section were to be carried to the continuous setting, it would mean that solutions to $\eqref{eq:hild}$ are bounded in the H\"older class $C^{1/2}$. This regularization effect has recently been conjectured in \cite{Vicol} for the vanishing viscosity approximation. If such a bound were to hold, it would mean that $\alpha=1/4$ is the true critical power for regularity and would answer the open question regarding $\eqref{eq:hild}$ stated in the introduction. For the SQG equation $\eqref{eq:SQG}$, the $C^{1-2\alpha}$ norm is critical \cite{DongP}: weak solutions that are bounded in time in $C^{1-2\alpha}$ are classical solutions.

\smallskip
First, we will prove bounds on $\eqref{eq:DE}$, the model without dissipation, as it is more elementary. Then we will generalize to the full model $\eqref{eq:DDE}$ with dissipation.

\bigskip
We rewrite $\eqref{eq:DE}$ in a different form, which allows us to give a more detailed picture
of how blow-up can occur. Let $b_k=(a_k-a_{k-1})2^k$. Then the $b_k$'s satisfy 
$$
b_k'(t)=-b_k^2+2b_{k-1}^2
$$
with corresponding initial data
$$
b_k(0)=-(a_k^0-a_{k-1}^0)2^k:=b_k^0
$$
By convention, we set $b_1'=-b_1^2$ and $b_0=0$.
We will take $b_k(0)>0$ for all $k$ so $b_k$ will remain pnon-negative by Lemma 2.4(a).
The fact that the $a_k$'s blow up in $X^s$, $s>1$, means that there exists $T>0$ such that $\lim_{t\to T} \sup_k b_k(t)2^{k(s-1)} =\infty$. In what follows, we work with $X^s$, $s>1$ solutions. We will need the following group of lemmas.

\bigskip
\noindent
\begin{lemma}
Suppose that for $T_1\le t\le T_2$, $b_{k-1}(t)<\sqrt{2}b_{k-2}(t)$ (or $>$) and $b_k(T_1)<\sqrt{2}b_{k-1} (T_1)$ (or $>$). Then for $T_1< t<T_2$, we have $b_k(t) < \sqrt{2}b_{k-1}(t)$ (or $>$).
\end{lemma} 

\medskip
\noindent
{\bf Proof} 

\medskip
For a contradiction suppose the set $A:= \{T_2\ge t > T_1: b_k(t)\ge \sqrt{2}b_{k-1}(t)\}$ is non-empty. By continuity, $A$ has an minimum call it $t_0$. Then $b_k(t_0)=\sqrt{2}b_{k-1}(t_0)$, which implies
\begin{eqnarray*}
(b_k-\sqrt{2}b_{k-1})'(t_0) &=& -b_k^2(t_0)+2b_{k-1}^2(t_0)+\sqrt{2}b_{k-1}^2(t_0)-2\sqrt{2}b_{k-2}^2(t_0) \\ &=& \sqrt{2}(b_{k-1}^2(t_0)-2b_{k-2}^2(t_0)) <0.
\end{eqnarray*}
From this, we can arrive at a contradiction to $t_0$ being the minimum of $A$. The proof 
for the opposite inequality is similar $\Box $

\medskip
The proof of the following lemma can be done in a similar way.
\bigskip
\noindent
\begin{lemma}
Suppose that for $T_1\le t\le T_2$, $b_{k-1}(t)< b_{k-2}(t)$ (or $> $) and $b_k(T_1)<b_{k-1}(T_1)$ (or $>$). Then for $T_1\le t<T_2$, we have $b_k(t) < b_{k-1}(t)$ (or $>$).
\end{lemma}

\begin{comment}
\medskip
\noindent
The following lemma is not used later, but may be of interest.
\bigskip
\noindent
\begin{lemma}
Let $1<q<\sqrt{2}$. Suppose that for $T_1\le t\le T_2$, $b_{k-1}(t)>q b_{k-2}(t)$ and $b_k(T_1)>qb_{k-1}(T_1)$, then $b_k(t) > qb_{k-1}(t)$ for $T_1<t<T_2$.
\end{lemma}
\end{comment}

\bigskip
Now, we analyze the following scenario. Suppose there exists $K_0$ such that the initial data $\{b_k^0\}$
satisfies
\begin{eqnarray}
\label{eq:1}
b_{k}^0 & < & b_{k-1}^0, \quad  \mbox{for}\,\, k> K_0 \\
\label{eq:2}
b_{k-1}^0<b_k^0 & <& \sqrt{2} b_{k-1}^0,\quad \mbox{for}\,\, k\le K_0
\end{eqnarray}
\begin{figure}[b!]
\includegraphics[width=0.6\linewidth]{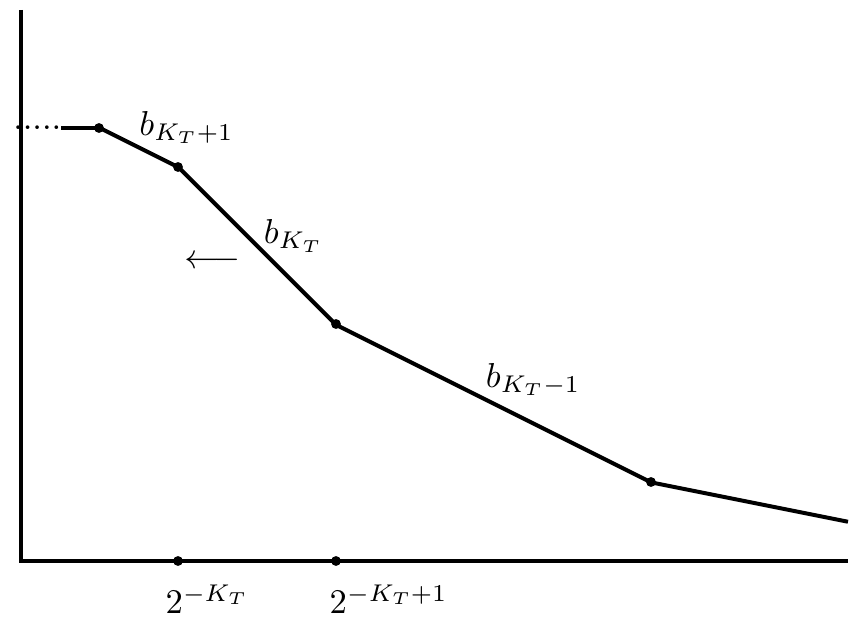}
\caption{}
\end{figure}
and consider solutions $\{b_k\}$ with this type of initial data. Observe that $\max_k b_k^0=b_{K_0}$. As we let time evolve, due to the control given by the  above lemmas, the solution will have a similar structure as the initial data. That is, at anytime $t<T$, there exists $K_t$ such that
\begin{eqnarray*}
b_{k}(t) & < & b_{k-1}(t), \quad  \mbox{for}\,\, k> K_t \\
b_{k-1}(t)<b_k(t) & <& \sqrt{2} b_{k-1}(t),\quad \mbox{for}\,\, k\le K_t.
\end{eqnarray*} 
At this point, we have $\max_{k} b_k(t)= b_{K_t}(t)$ and $K_t\ge K_0$. Thinking in terms of $a_k\approx \theta(2^{-k},t)$ and $b_k$ being the modulus of the slope of the line between $(2^{-k},a_k)$ and $(2^{-k+1},a_{k-1})$, we are approximating a solution $\theta$ of $\eqref{eq:hild}$ by a piecewise linear function taking the value $a_k$ at $2^{-k}$ (see figure 2).

\begin{comment}
\begin{figure}[b!]
\includegraphics[width=0.6\linewidth]{picture.pdf}
\caption{}
\end{figure}
\end{comment}

\bigskip
\noindent
The picture is of a front propagating towards the origin attempting to form a cusp. Once the front hits the origin (i.e. the time when $\sup_k b_k=\lim_{k\to \infty} b_k$), blow-up occurs. The index $K_t$ can be thought of as the position of the inflection point of the function that sequence $\{b_k\}$ is approximating. One can hypothesize that $b_k(t) \approx \sqrt{2}b_{k-1}(t)$ as $t\to T$. In fact, the following result suggests that this may be the case.

\begin{theorem}
Suppose we have initial data satisfying $\eqref{eq:1}$ and $\eqref{eq:2}$. Then for $0<\beta\le1/2$,
$$
\sup_{t<T, k} b_k2^{k(\beta-1)} <\infty
$$
\end{theorem}

\Pf Follows from Lemma 3.1 and 3.2. $\Box $

\begin{comment}
% wrong
Observe that for $t<T$,
$$
b_{K_0}'(t)=-b_{K_0}^2(t)+2b_{K_0-1}^2(t) \le -b_{K_0}^2(t)+2(b_{K_0-1}^0)^2.
$$
Then by Lemma 0.3 applied with $c=2(b_{K_0-1}^0)^2$ , $b_{K_0}(t) < \sqrt{2}b_{K_0-1}^0$. Now,
$$
b_{K_0+1}'=-b_{K_0+1}^2+2b_{K_0}^2 \le -b_{K_0+1}^2+2(\sqrt{2}b_{K_0-1}^0)^2
$$
so by Lemma 0.3, $b_{K_0+1}(t)<2b_{K_0-1}^0$.
From an induction argument, it is not hard to
see that for $k>K_0$, we have $b_k(t) < (\sqrt{2})^{k-K_0+1}b_{K_0-1}^0$. The result follows. $\Box$
\end{comment}

\subsection{With Dissipation}

\medskip

We will consider a class of initial data slightly different from the previous section. The blow-up picture becomes less clear, but we still are able to have the analogue to Theorem 3.3. We still consider initial data $\{a_k^0\}_{k=0}^\infty\in X^s$, $s>1$, that is positive and non-decreasing in $k$. As before we let $b_k=(a_k-a_{k-1})2^k$ and we start with the following configuration
\begin{eqnarray}
b_{k}^0 & < & b_{k-1}^0, \quad  \mbox{for}\,\, k> K_0 \\
b_{k}^0 & <& \sqrt{2}b_{k-1}^0,\quad \mbox{for}\,\, k\le K_0 \\
\sup_k b_k^0 &=& b_{K_0}^0
\end{eqnarray}
for some positive integer $K_0$. The following lemma is an analogue to the observation made with the model without dissipation.

\begin{lemma}
Let $0<\alpha<1/2$. Suppose we have initial data as in (11)-(13). For all $t$ such that the solution is smooth, there exists a positive integer $K_t$ such that
\begin{eqnarray}
b_{k} & < & b_{k-1}, \quad\qquad  \mbox{for}\,\, k> K_t \\
b_{k} & <& \sqrt{2}b_{k-1},\quad \mbox{for}\,\, k\le K_t 
\end{eqnarray}
In addition, $K_t \ge K_0$.
\end{lemma}

\Pf The result will quickly follow from the following two claims.

\bigskip
{\bf Claim 1}: Fix $t$. Suppose we have a $j$ such that $b_j(t)<b_{j-1}(t)<b_{j-2}(t)$. For later times $r>t$, as long as $b_{j-1}(r)<b_{j-2}(r)$, we have $b_j(r)<b_{j-1}(r)$.

\medskip
{\bf Claim 2}: Fix $t$. Suppose we have a $j$ such that $b_j(t)<\sqrt{2}b_{j-1}(t)$ and $b_{j-1}(t)<\sqrt{2}b_{j-2}(t)$. For later times $r>t$, as long as $b_{j-1}(r)<\sqrt{2} b_{j-2}(r)$, we have $b_j(r)<\sqrt{2}b_{j-1}(r)$.

\medskip
\ul{Proof of Claim 1:}

\medskip
For a contradiction, suppose there exists a minimal time $t_0>t$ such that $b_j(t_0)=b_{j-1}(t_0)$ or equivalently $(a_j-a_{j-1})(t_0)=\frac{1}{2} (a_{j-1}-a_{j-2})(t_0)$. First,
\begin{eqnarray*}
(b_j-b_{j-1})'(t_0) &=& -b_j^2+b_{j-1}^2 +2(b_{j-1}^2-b_{j-2}^2) + 2^j\left( (\mathscr{L} a)_{j-1}- (\mathscr{L} a)_j\right) + 2^{j-1}\left( (\mathscr{L} a)_{j-2}- (\mathscr{L} a)_{j-1}\right) \\
&=& 2(b_{j-1}^2-b_{j-2}^2) - 2^j\left( (\mathscr{L} a)_{j}- (\mathscr{L} a)_{j-1}\right) + 2^{j-1}\left( (\mathscr{L} a)_{j-1}- (\mathscr{L} a)_{j-2}\right)
\end{eqnarray*}
where everything from now on will be evaluated at $t=t_0$.
By hypothesis, the first term on the right side is negative. If we can show
$$
(\mathscr{L} a)_j- (\mathscr{L} a)_{j-1} > \frac{1}{2} \left( (\mathscr{L} a)_{j-1}- (\mathscr{L} a)_{j-2}\right)
$$
we will have a contradiction. From now on, we are left with lengthy algebraic computations. By direct calculation,
\begin{eqnarray}
\label{eq:dissdiff}
(\mathscr{L} a)_j - (\mathscr{L} a)_{j-1} &=& \sum_{n=0}^{j-2} (a_j-a_{j-1})2^{2\alpha n} + 2(a_j-a_{j-1})2^{2\alpha (j-1)}  \\ \nonumber &\, & + 2^{(1+2\alpha)j} (1-2^{-(1+2\alpha)})\sum_{n=j+1}^\infty  (a_j-a_n)2^{-n}
\end{eqnarray}
Also, we have
\begin{eqnarray*}
\frac{1}{2}\left[(\mathscr{L} a)_{j-1}- (\mathscr{L} a)_{j-2}\right] &=& \sum_{n=0}^{j-3} \frac{1}{2}(a_{j-1}-a_{j-2})2^{2\alpha n} + (a_{j-1}-a_{j-2})2^{2\alpha (j-2)}  \\ &\, & + 2^{(1+2\alpha)(j-1)-1}(1-2^{-(1+2\alpha)}) \sum_{n=j}^\infty  (a_{j-1}-a_n)2^{-n} \\
&=& \sum_{n=0}^{j-3} (a_{j}-a_{j-1})2^{2\alpha n} + 2 (a_{j}-a_{j-1})2^{2\alpha (j-2)}  \\ &\, & + 2^{(1+2\alpha)(j-1)-1}(1-2^{-(1+2\alpha)}) \sum_{n=j}^\infty  (a_{j-1}-a_n)2^{-n}
\end{eqnarray*}
where we have used $b_k=b_{k-1}$ in the last equality. Then
\begin{eqnarray*}
(\mathscr{L} a)_j - (\mathscr{L} a)_{j-1}  - \frac{1}{2}\left[(\mathscr{L} a)_{j-1}- (\mathscr{L} a)_{j-2}\right] &=& (2-2^{-2\alpha})(a_j-a_{j-1})2^{2\alpha(j-1)} \\ &\, & + \frac{1}{4}(1-2^{-(1+2\alpha)})(a_j-a_{j-1}) 2^{2\alpha(j-1)} \\ &\,& +2^{(1+2\alpha)j} (1-2^{-(1+2\alpha)})\sum_{n=j+1}^\infty  (a_j-a_n)2^{-n} \\
&\,& -2^{(1+2\alpha)(j-1)-1}(1-2^{-(1+2\alpha)}) \sum_{n=j+1}^\infty  (a_{j-1}-a_n)2^{-n} \\
&=& (2+2^{-2} -2^{-2\alpha} - 2^{-(3+2\alpha)})(a_j-a_{j-1}) 2^{2\alpha(j-1)} \\
&\, & +(1-2^{-(1+2\alpha)}) 2^{(1+2\alpha)(j-1)-1} \sum_{n=j+1}^\infty (a_j-a_{j-1})2^{-n} \\
&\, & +(1-2^{-(1+2\alpha)})(1- 2^{-(2+2\alpha)}) 2^{(1+2\alpha) j} \sum_{n=j+1}^\infty (a_j-a_n)2^{-n}.
\end{eqnarray*}
Using that $b_{k+1}< b_k$ for $k\ge j$, we have that $\ds \sum_{n=j+1}^\infty (a_j-a_n)2^{-n} \ge -\frac{2}{3}(a_j-a_{j-1})2^{-j}$. Inserting this into the above estimate we get that the above expression is bounded below by
$$
\left( \frac{5}{2} -5(2^{-(2+2\alpha)}) -\frac{2^{1+2\alpha}}{3} (1-2^{-(1+2\alpha)})(1-2^{-(2+2\alpha)})\right)(a_j-a_{j-1} )2^{2\alpha(j-1)}.
$$
Since $0<\alpha<1$, one can see that the expression is indeed positive.

\medskip
\ul{Proof of Claim 2:} The proof is done similarly to Claim 1, but we will provide the details for completeness. Assume there exists a time $t_0$ such that $b_j(t_0)=\sqrt{2} b_{j-1}(t_0)$. Also, at this time, $b_k(t_0) \le \sqrt{2}b_{k-1}(t_0)$ for all $k\ne j$. If we can show $(\sqrt{2}b_{j-1}-b_j)'(t_0)>0$ then we have a contradiction and the claim is proven. From now on all function will be evaluated at $t_0$. By direct calculation,
\begin{eqnarray*}
(\sqrt{2}b_{j-1}-b_j)' &=& (b_j^2 -\sqrt{2} b_{j-1}^2) - 2(b_{j-1}^2-\sqrt{2}b_{j-2}^2) \\
&\,& + [((\mathscr{L} a)_j-(\mathscr{L} a)_{j-1})- \frac{1}{\sqrt{2}} ((\mathscr{L} a)_{j-1}-(\mathscr{L} a)_{j-2})] 2^j
\end{eqnarray*}
The contribution on the first line is positive. It suffices to show 
$$
D_j:=((\mathscr{L} a)_j-(\mathscr{L} a)_{j-1})- \frac{1}{\sqrt{2}} ((\mathscr{L} a)_{j-1}-(\mathscr{L} a)_{j-2})>0
$$
After similar algebraic manipulations, using that $b_j=\sqrt{2}b_{j-1}$ as in Claim 1 we get 
$$
D_j =(1-2^{-(1+2\alpha)})\left[ \left(\frac{2\sqrt{2}+1}{\sqrt{2}}\right)(a_j-a_{j-1})2^{2\alpha(j-1)} + (1-2^{-(3/2+2\alpha)})2^{(1+2\alpha)j} \sum_{n=j+1}^\infty (a_j-a_n)2^{-n}\right]
$$
Using $b_k \le \sqrt{2} b_{k-1}$ for $k\ne j$, 
$$
 \sum_{n=j+1}^\infty (a_j-a_n)2^{-n} \ge  - \frac{2}{2\sqrt{2}-1} (a_j-a_{j-1})2^{-j}
 $$
 Inserting this bound into the expression for $D_j$, we see that $D_j$ is indeed positive for $0<\alpha<1/2$.

\medskip
Thus, by Claims $1$ and $2$, we know the solution must be of the form given by (13) and (14). At time $t=0$, we have control over all $b_k$'s except for $b_{K_0+1}$. Though we initially have $b_{K_0+1}(0)< b_{K_0}(0)$, there may exist a time $t_1$ such that $b_{K_0}(t_1)\le b_{K_0+1}(t_1) < \sqrt{2}b_{K_0}(t_1)$. By Claim 2, we know for $t>t_1$, $b_{K_0+1}(t) < \sqrt{2}b_{K_0}(t)$. In addition, by Claim 1, at time $t_1$, our solution is of the form (13) and (14) with $K_{t_1} = K_0+1$. If there exists no time $t_1$, then $K_t=K_0$. Continuing with this line of reasoning, the result follows. $\Box$

\begin{theorem}
(A-priori H\"older-$1/2$ bound) Let $\{a_k\}$ be a solution to $\eqref{eq:DDE}$ with initial data $\{a_k^0\}\in X^s$, $s>1$, that is increasing in $k$ and positive with $b_j^0< \sqrt{2}b_{j+1}^0$ for $j\ge 0$. Then 
$$
\sup_{t>0,k} (a_k(t)-a_{k-1}(t))2^{k/2} = \sup_{t>0, k} b_k(t) 2^{-k/2} < \infty
$$

\end{theorem}

\Pf By Claim 2 in the proof of the previous Lemma, $b_k(t)< \sqrt{2}b_{k+1} (t)$ for all $t$ and $k\ge 0$. Then for all $k$, $b_k(t) < 2^{\frac{k-1}{2}} b_1(t)< C 2^{k/2}$. $\Box$

\section{Regularity for $1/4< \alpha \le 1/2$ and Blow-up for $\alpha<1/4$}

\medskip
\noindent
The previous theorem will allow us to adapt the argument of Theorem 2.5 and prove blow-up for $\eqref{eq:DDE}$ for $0<\alpha<1/4$.

\bigskip
\begin{theorem} Let $0<\alpha<1/4$ Suppose $\{a_k^0\}$ is increasing in $k$ and is non-negative. Then there exists initial datum in $X^s$, $s>1$, such that solutions blow-up in finite time.
\end{theorem}

\Pf

\medskip
\noindent
As before we let 
$$
J(t) = \sum_{k=0}^\infty (a_\infty(t)-a_k(t)) 2^{k\delta}.
$$
for $0<\delta<1-4\alpha$ where $a_\infty(t)=\lim_{k\to\infty} a_k(t)$. 
Then computing and using Lemma 2.6 we get
\begin{eqnarray*} 
\frac{d}{dt} J(t) &=& \sum_{k=1}^\infty (a_k(t)-a_{k-1}(t))^2 2^{k(\delta+1)} - \sum_{k=1}^\infty \left( (\mathscr{L} a)_\infty - (\mathscr{L}a)_k\right)2^{k\delta}\\
& \ge & C_0 \sum_{k=1}^\infty (a_\infty(t) -a_k(t))^2 2^{k(\delta+1)} - \sum_{k=1}^\infty \left( (\mathscr{L} a)_\infty - (\mathscr{L}a)_k\right)2^{k\delta}
\end{eqnarray*}
where 
$$
(\mathscr{L} a)_\infty = \lim_{k\to \infty} (\mathscr{L}a)_k = \sum_{n=0}^\infty (a_\infty-a_n) 2^{2\alpha n}
$$
Now,
\begin{eqnarray*}
(\mathscr{L} a)_\infty - (\mathscr{L}a)_k = \sum_{n=0}^{k-1}(a_\infty -a_k)2^{2\alpha n} + \sum_{n=k}^\infty (a_\infty -a_n)2^{2\alpha n} - \sum_{n=k+1}^\infty (a_k-a_n)2^{(1+2\alpha)k} 2^{-n}
\end{eqnarray*} 
We analyze each of the above three terms separately.
$$
\sum_{n=0}^{k-1}(a_\infty -a_k)2^{2\alpha n} \le C (a_\infty - a_k) 2^{2\alpha k}
$$
By Theorem 3.5,
$$
 \sum_{n=k}^\infty (a_\infty -a_n)2^{2\alpha n} \le C \sum_{n=k}^\infty 2^{(2\alpha -1/2)n} \le C2^{(2\alpha-1/2)(k-1)}
 $$
 For the last term we have
 $$
 \sum_{n=k+1}^\infty (a_n-a_k)2^{(1+2\alpha)k} 2^{-n} \le (a_\infty -a_k) 2^{(1+2\alpha)k}\sum_{n=k+1}^\infty 2^{-n} \le C (a_\infty -a_k)2^{2\alpha k}.
 $$
 Then using $\alpha<1/4$ and $0<\delta<1-4\alpha$,
 $$
 \sum_{k=1}^\infty \left( (\mathscr{L} a)_\infty - (\mathscr{L}a)_k\right)2^{k\delta} \le C_1+C_2\sum_{k=1}^\infty (a_\infty -a_k) 2^{(2\alpha+\delta)k}
 $$
 For every $\epsilon>0$, by an application of Holder's inequality and using that $\alpha <1/4$,
 $$
 \sum_{k=1}^\infty (a_\infty -a_k)2^{(2\alpha+\delta)k} \le \epsilon\sum_{k=1}^\infty (a_\infty -a_k)^2 2^{(\delta+1)k} + C_\epsilon \|a\|_{\ell^\infty}
 $$
Given these bounds and following the proof of Theorem 2.6 we can show
 $$
 J'(t) \ge C(\delta) J(t)^2 - C(1+ \|a\|_{\ell^\infty}).
 $$
By choosing initial data appropriately, the inequality above leads to blow-up. $\Box$

\medskip
Now, we will move towards proving regularity for $1/4<\alpha<1/2$.

\begin{theorem} Let $1/4<\alpha<1/2$. Suppose we have the same hypotheses as in Theorem 3.5. Then solutions exist for all time, in particular, for $s> 1$
$$
\sup_{t>0} \|a(t)\|_{X^s} < \infty.
$$
\end{theorem}

\Pf 

\medskip
After a computation, 
\begin{equation}
\label{eq:bound}
(b_{k,s})'(t) = -b_k\cdot b_{k,s}(t)+2^s b_{k-1}\cdot b_{k-1,s} (t)- \left((\mathscr{L} a)_k-(\mathscr{L}a)_{k-1}\right)(t)2^{sk}
\end{equation}
From Theorem 3.5, $b_{k-1}(t) \le C_1 2^{k/2}$ for all $t\in [t_0,t_1]$ where the constant $C_1$ is independent of the time interval. Then there exists $K'$ such that if $k\ge K'$ and $b_{k,s}(t)> c_s \|a\|_{X^s}(t)$ where $c_s$ is the constant from Lemma 2.1, 
$$
(b_{k,s})'(t)  \le b_{k,s}\left(C_1 2^{k/2+s} - C(\alpha)(2^{2\alpha k}-2^{2\alpha})\right)<0.
$$
In the above inequality, we use $0<\alpha<1/4$. This implies that for any $T>0$, $\lim_{t\to T, k\to\infty} b_{k,s}(t)\ne \infty$. Therefore, the solution exists globally in time. $\Box$

\bigskip
\noindent
{\bf Remark:} It is unclear whether the hypothesis that $b_j^0 \le \sqrt{2}b_{j+1}^0$ can be weakened as it is crucial to the proof of the a-priori bound. We conjecture that any reasonable initial data will eventually satisfy such a condition at least for large enough indexes.

\bigskip
\noindent
{\bf Summary}: We have shown that when $\alpha>1/4$, under mild assumptions on the initial data, it is possible to have global solutions: the dissipation terms win over the nonlinearity. In order to have blow-up, energy initially present in the lower modes must reach the higher modes. By making use of a-priori Holder-type control and an estimate on the dissipative terms, we have shown that such an energy transfer must stop at some time. 

\medskip
It is important to note that our conditions on the initial data are not a condition of ``smallness'' with regard to some norm. Our results of global regularity can hold for large initial data. Also, one can easily find initial data for which our model blows up for $0<\alpha<1/4$ and is regular for $\alpha>1/4$.

\section{Acknowledgements}
This work was completed under the guidance of Prof. Alex Kiselev and with support from NSF-DMS grants 1412023 and 1453199. The author also wishes to thank Vlad Vicol and Xiaoqian Xu for useful discussions.

\bibliographystyle{plain}
\bibliography{discreterefs}

\end{document}